%

\magnification=\magstep1
\def\forces{\parallel\!\!\! -}
\def\restrict{{\restriction}}
\def\Smallskip{\vskip1.5truecm}
\def\Bigskip{\vskip2.5truecm}

\def\qed{{\vcenter{\hrule height.4pt \hbox{\vrule width.4pt height5pt
 \kern5pt \vrule width.4pt} \hrule height.4pt}}}
\def\ok{\vbox{\hrule height 8pt width 8pt depth -7.4pt
    \hbox{\vrule width 0.6pt height 7.4pt \kern 7.4pt \vrule width 0.6pt height 7.4pt}
    \hrule height 0.6pt width 8pt}}
\def\nt{{\leq}\kern-1.5pt \vrule height 6.5pt width.8pt depth-0.5pt \kern 1pt}
\def\sd{{\times}\kern-2pt \vrule height 5pt width.6pt depth0pt \kern1pt}
\def\notin{{\in}\kern-5.5pt / \kern1pt}
\def\ZZ{{\rm Z}\kern-3.8pt {\rm Z} \kern2pt}
\def\RR{{\rm I\kern-1.6pt {\rm R}}}
\def\NN{{\rm I\kern-1.6pt {\rm N}}}
\def\QQ{{\Bbb Q}}
\def\CC{{\Bbb C}}
\def\DD{{\Bbb D}}
\def\AA{{\Bbb A}}
\def\BB{{\Bbb B}}
\def\MM{{\Bbb M}}
\def\PP{{\Bbb P}}
\def\KK{{\rm I\kern-1.6pt {\rm K}}}
\def\11{{\rm 1}\kern-2.2pt {\rm \vrule height6.1pt
    width.3pt depth0pt} \kern5.5pt}
\def\zp#1{{\hochss Y}\kern-3pt$_{#1}$\kern-1pt}

\def\egs{\vrule height 6pt width.5pt depth 2.5pt \kern1pt}
\font\capit=cmcsc10 scaled\magstep0
\font\capitg=cmcsc10 scaled\magstep1
\font\dunh=cmdunh10 scaled\magstep0
\font\dunhg=cmr10 scaled\magstep1
\font\dunhgg=cmr10 scaled\magstep2

\font\sanse=cmss10 scaled\magstep0
\overfullrule=0pt
\openup1.5\jot
\input mssymb

\centerline{}
\Bigskip
\centerline{\dunhgg Perfect sets of random reals}
\Bigskip
\centerline{\capitg J\"org Brendle\footnote{$^*$}{{\rm
The first author would like to thank the MINERVA-foundation
for supporting him}} and Haim Judah\footnote{$^{**}$}{{\rm
The second author would like to thank the Basic Research
Foundation (the Israel Academy of Sciences and Humanities) for
supporting him}}}
\Smallskip
\centerline{Abraham Fraenkel Center for Mathematical Logic} \par
\centerline{Department of Mathematics} \par
\centerline{Bar--Ilan University} \par
\centerline{52900 Ramat--Gan, Israel} 
\Bigskip
\centerline{\capit Abstract}
\bigskip
\noindent We show that the existence of a perfect set of random
reals over a model $M$ of $ZFC$ does not imply the existence of
a dominating real over $M$, thus answering a well-known open question 
(see [BJ 1] and [JS 2]). We also prove that $\BB \times 
\BB$ (the product of two copies of the random algebra) neither
adds a dominating real nor adds a perfect set of random reals (this
answers a question that A. Miller asked during the logic
year at MSRI).
\vfill\eject
\noindent{\dunhg Introduction}
\Smallskip
The goal of this work is to give several results concerning
the relationship between perfect sets of random reals,
dominating reals, and the product of two copies of the
random algebra $\BB$. Recall that $\BB$ is the algebra
of Borel sets of $2^\omega$ modulo the null sets. Also,
given two models $M \subseteq N$ of $ZFC$, we say that
$g \in \omega^\omega \cap N$ is a {\it dominating real over
$M$} iff $\forall f \in \omega^\omega \cap M \;
\exists m \in \omega \; \forall n \geq m \; (g(n) > f(n))$;
and $r \in 2^\omega \cap N$ is {\it random over $M$}
iff $r$ avoids all Borel null sets coded in $M$ iff $r$
is the real determined by some filter which is $\BB$-generic
over $M$ (see [Je 1, section 42] for details).
\par
A tree $T \subseteq 2^{< \omega}$ is {\it perfect} iff $\forall t 
\in T \; \exists s \supseteq t \; (s \hat{\;} \langle 0
\rangle \in T \; \land s \hat{\;} \langle 1 \rangle \in T)$.
For a perfect tree $T$ we let $[T] := \{ f \in 2^\omega ;
\; \forall n \; (f \restrict n \in T ) \}$ denote the
set of its branches. Then $[T]$ is a perfect set (in the
topology of $2^\omega$). Conversely, given a perfect
set $S \subseteq 2^\omega$ there is perfect tree $T \subseteq
2^{< \omega}$ such that $[T] = S$. This allows us to confuse
perfect sets and perfect trees in the sequel; in particular,
we shall use the symbol $T$ for both the tree and the set of its branches.
--- As a perfect tree is (essentially) a real, the statement
{\it there is a perfect set of reals random over $M$ in
$N$} (where $M \subseteq N$ are again models of $ZFC$) asserts
the existence of a certain kind of real in $N$ over $M$;
and thus we may ask how it is related to the existence
of other kinds of reals (like dominating reals). This will be our main topic. ---
We recall that the existence of a random real does not 
imply the existence of a perfect set of random reals; in
fact Cicho\'n showed that $\BB$ does not add a perfect
set of random reals [BJ 1, Theorem 2.1]. (Here, we say that
a p.o. $\PP$ {\it adds a perfect set of random reals} iff
there is a perfect set of reals random over $M$ in $M[G]$,
where $G$ is $\PP$-generic over $M$; a similar definition
applies to dominating reals etc.)
\par
We note that being a perfect set of random reals over
some model $M$ of $ZFC$ is absolute in the following sense: if
$M \subseteq N_0 \subseteq N_1$ are models of $ZFC$,
$T \in (2^{< \omega})^\omega \cap N_0$ is a perfect tree so that
$[T] \cap {N_0}$ consists only of reals random over $M$, then
every real in $[T] \cap {N_1}$ is random over $M$ as well
(see [Je 1, Lemma 42.3]).
\par
We now state the main results of our work, and explain how
they will be presented in $\S\S$ 1 -- 3; then we will give some further
motivation for the study of perfect sets of random reals, and
close with some notation.
\bigskip
{\sanse The main results and the organization of the paper.}
Using techniques of [Ba], Bartoszy\'nski and Judah proved
in [BJ 1, Theorem 2.7] that
\smallskip
\item{(*)} {\it given models of $ZFC$ $M \subseteq N$ such that
$N$ contains a dominating real over $M$, $N[r]$ contains
a perfect set of random reals over $M$, where $r$ is random
over $N$.}
\smallskip
\noindent Our first result shows that the converse does not hold
(1.4 -- 1.7).
\smallskip
{\capit Theorem 1.} {\it There is a p.o. which adds a perfect
set of random reals and does not add dominating reals.}
\smallskip
\noindent The framework for proving Theorem 1 (developed in
$\S$ 1) will enable us to give general preservation results
for {\it not adding dominating reals} in both finite
support iterations and finite support products
of $ccc$ forcing notions (1.8). The former will be exploited
in 1.9 to discuss cardinal invariants closely related
to our subject. As a special instance of the latter
we shall show (1.10)
\smallskip
{\capit Theorem 2.} {\it $\BB \times \BB$ does not add dominating
reals.}
\smallskip
\noindent (This result was proved earlier by Shelah but never
published.) The algebra $\BB \times \BB$ is rather different
from $\BB$; e.g. it is well-known that $\BB \times \BB$ adds
Cohen reals whereas $\BB$ does not (see [Je 2, part I, 5.9]).
As it is known that some other forcing notions adding both
Cohen and random reals (like $\BB \times \CC \cong \BB *
\CC$ and $\CC * \BB$, where $\CC$ is the Cohen algebra, and
$*$ denotes iteration) do not add perfect sets of random
reals (see [JS 2, 2.3] for $\BB \times \CC$ and [BJ 1, Theorem
2.13] for $\CC * \BB$), we may ask whether $\BB \times \BB$
does. We shall show in $\S$ 2 that the answer is again negative.
\smallskip
{\capit Theorem 3.} {\it $\BB \times \BB$ does not add a perfect set
of random reals.}
\smallskip
\noindent The argument for this proof (which uses ideas from
the proof that $\BB$ does not add a perfect set of random
reals --- see [BJ 1, 2.1 -- 2.4]) is rather long and technical;
and one might get a shorter proof if the following
question has a positive answer.
\smallskip
{\dunh Question 1.} {\it Is $\BB \times \BB$ a complete subalgebra
of $\CC * \BB$?}
\smallskip
\noindent We note here that all other embeddability relations
between these three algebras adding both Cohen and random
reals (namely, $\BB \times \CC$, $\BB \times \BB$,
$\CC * \BB$) are known. We shall sketch the arguments which
cannot be found in literature in 3.1.
\par
Two further open problems are closely tied up with the
Bartoszy\'nski--Judah Theorem (*) and our Theorem 1,
respectively. 
\smallskip
{\dunh Question 2.} {\it Given models of $ZFC$ $M \subseteq N$
such that $N$ contains both a dominating real and a random
real over $M$, is there a perfect set of random reals over
$M$ in $N$?}
\smallskip
\noindent We shall show in 3.2 that to answer Question 2
it suffices to consider the problem whether $\BB \times
\DD \cong \BB * \DD$ adds a perfect set of random reals, where
$\DD$ is Hechler forcing. We note that for many
p.o.s $\PP$ adding a dominating real 
(e.g. Mathias forcing) it is true that $\BB
\times \PP$ adds a perfect set of random reals (3.3).
\smallskip
{\dunh Question 3.} {\it Given models of $ZFC$ $M \subseteq N$,
does the existence of a perfect set of random reals over
$M$ in $N$ imply the existence of an unbounded real over 
$M$ in $N$?}
\smallskip
\noindent Here we say that $g \in \omega^\omega \cap N$ is
an {\it unbounded real over} $M$ iff $\forall f \in
\omega^\omega \cap M \; \exists^\infty n \; (g(n) > f(n))$,
where $\exists^\infty n$ means {\it there are infinitely
many $n$} (dually, $\forall^\infty n$ abbreviates
{\it for all but finitely many $n$}).
\bigskip
{\sanse Motivation.} One of the reasons for studying perfect
sets of random reals concerns finite support iterations
of $ccc$ forcing notions. Namely, let $\langle \PP_n ,
\breve Q_n ; \; n \in \omega \rangle$ be an $\omega$-stage
finite support iteration such that for all $n \in \omega$,
$\forces_{\PP_n} $ "$\breve Q_n$ is $ccc$". Then the 
following are equivalent [JS 2, Theorem 2.1]:
\smallskip
{\it \item{(i)} There exists $r \in V [G_\omega] \setminus
\bigcup_n V[G_n]$ random over $V$,
\par
\item{(ii)} there exists $n \in \omega$ and $T \in V[G_n]$
a perfect set of random reals over $V$, \par}
\smallskip
\noindent where $\langle G_i ; \; i \leq \omega \rangle$ is
a chain of $\PP_i$-generic filters. So {\it adding a random
real in the $\omega$-th stage} is stronger than just
{\it adding a random real in an initial step} (on the other
hand, {\it $\PP_\omega$ adds a dominating real} is simply
equivalent to {\it there is an $n \in \omega$ such
that $\PP_n$ adds a dominating real} [JS 1, Theorem 2.2]).
--- Also perfect sets of random reals seem to play
an important role in the investigation of the problem,
posed by Fremlin, whether the smallest covering of the
real line by measure zero sets can have cofinality
$\omega$. To build a model of $ZFC$ where this is true
we suggest an iterated forcing construction (with
finite supports) which firstly adds $\omega_\omega$
many Cohen reals over $L$ to produce a family of
$\omega_\omega$ null sets which will still cover
the real line in the final extension, and then goes through every
subalgebra of the random algebra which is the random
algebra restricted to some small inner model (in which
the continuum has size $< \omega_\omega$)
in $\omega_{\omega + 1}$ steps (see the introduction
of [JS 3] for details). By construction, 
we destroy all small covering families.
So the main problem is to show that we do not
add a real which does not belong to the family of $\omega_\omega$
null sets added in the intermediate stage. To do this, it suffices
(essentially) to prove that the whole iteration does not add a
perfect set of random reals over the ground model $L$.
We think that our Theorem 3 is a small but important
step in this direction, and we hope that the ideas
involved can be generalized to give a positive answer to
\smallskip
{\dunh Question 4.} {\it Let $\AA$ be a complete subalgebra
of $\BB$. Let $\breve \BB_\AA$ be an $\AA$-name for $\BB$.
Is it true that $\BB * \breve \BB_\AA$ does not add a
perfect set of random reals?}
\smallskip
\noindent Cicho\'n's Theorem [BJ 1, Theorem 2.1] says that
this is true if $\AA = \BB$, and our Theorem 3 gives a positive
answer in case $\AA$ is trivial.
\bigskip
{\sanse Notation.} Our notation is fairly standard. We refer
the reader to [Je 1] for set theory and to [Ox] for
measure theory. Most of the cited material will appear in
the forthcoming book [BJ 2]. We now explain some notions
which might be less familiar. \par
Given a finite sequence $s$ (i.e. either
$s \in 2^{<\omega}$ or $s \in \omega^{<\omega}$), we let $lh(s) 
:= dom (s) $ denote the length of $s$; for $\ell \in lh(s)$,
$s \restrict \ell$ is the restriction of $s$ to $\ell$. 
$\hat{\;}$ is used for concatenation of sequences; and
$\langle \rangle$ is the empty sequence. Furthermore,
for $s \in 2^{< \omega}$, $[s] := \{ f \in 2^\omega ;
\; f \restrict lh (s) = s \}$ is the set of branches through
$s$ (the open subset of $2^\omega$ determined by $s$).
\par
Given a perfect tree $T \subseteq 2^{< \omega}$ and $s \in T$,
we let $T_s := \{ t \in T ; \; t \subseteq s$ or $s \subseteq t
\}$; and $stem(T) := \cup \{ s \in T ; \; T_s =T \}$ is the
stem of $T$. For $\ell \in \omega$, we let $T \restrict \ell
:= \{ s \in T ; \; lh (s) \leq \ell \}$, the finite initial
part of $T$ of height $\ell$. We will confuse finite trees $T$
with all branches of fixed length $\ell$ with the set of
branches $[T] := \{ s \in T ; \; lh (s) = \ell \}$.
\par
We assume the reader to be familiar with forcing and Boolean-valued
models (see [Je 1], [Je 2]). We suppose that all our p.o.s
(forcing notions)
have a largest element $\11$. Given a p.o. $\PP \in V$, we
shall denote $\PP$-names by symbols like $\breve f$, $\breve T$, ...
and their interpretation in $V[G]$ (where $G$ is $\PP$-generic
over $V$) by $\breve f [G]$, $\breve T [G]$, ... If $\phi$ is a
sentence of the $\PP$-forcing language, we let $\parallel \phi
\parallel$ be the Boolean value of $\phi$; i.e. the maximal
element forcing $\phi$ in the complete Boolean algebra $r.o.(\PP)$
associated with $\PP$. We shall often confuse $\PP$ and $r.o.(\PP)$.
\par
We equip $\BB \times \BB : = \{ (p,q) ; \; p, q \in \BB \setminus
\{ 0 \} \} \cup \{ 0 \}$ with the product measure (i.e.
$\mu(p,q) = \mu (p) \cdot \mu (q)$). Then $\mu : \BB \times
\BB \to [0,1]$ is finitely additive and strictly positive (any
non-zero condition has positive measure). By [Ka, Proposition
2.1], $\mu$ can be extended to a finitely additive, strictly
positive measure on $r.o.(\BB \times \BB)$. This will be used
in 2.5. Note that this measure is not $\sigma$-additive.

\vfill\eject

\noindent{\dunhg $\S$ 1 Not adding dominating reals}
\Smallskip
{\sanse 1.1} We shall now introduce the framework needed to prove
Theorem 1. Besides giving the latter result this framework
will also provide us with preservation results for {\it not
adding dominating reals} in finite support products and finite
support iterations.
\par
Let $\PP$ be an arbitrary p.o. A function $h: \PP \to \omega$
is a {\it height function} iff $p \leq q$ implies $h(p) \geq h(q)$.
A pair $(\PP , h)$ is {\it soft} iff $\PP$ is a p.o., $h$ is
a height function on $\PP$, and the following two conditions
are satisfied: \par
\item{(I)} (decreasing chain property) if $\{ p_n ; \; n \in 
\omega \}$ is decreasing and $\exists m \in \omega \; \forall
n \in \omega \; (h(p_n) \leq m)$, then $\exists p \in \PP \;
\forall n \in \omega \; (p \leq p_n)$;
\par
\item{(II)} (weak finite cover property) given $m \in \omega$
and $\{ p_i ; \; i \in n \} \subseteq \PP$ there is $\{ q_j
; \; j \in k \} \subseteq \PP$ so that \par
\itemitem{(i)} $\forall i \in n , \; j \in k$, $q_j$ is
incompatible with $p_i$; \par
\itemitem{(ii)} whenever $q$ is incompatible with all $p_i$ and
$h(q) \leq m$ then there exists $j \in k$ so that $q \leq q_j$.
\par
\noindent We also consider the following property of pairs
$(\PP , h)$ --- where $\PP$ is a p.o. and $h$ a height function
on $\PP$: \par
\item{(*)} given a maximal antichain $\{ p_n ; \; n \in \omega \}
\subseteq \PP$ and $m \in \omega$ there exists $n \in \omega$
such that: whenever $p$ is incompatible with $\{ p_j ; \;
j \in n \}$ then $h(p) > m$. \bigskip
{\sanse 1.2} {\capit Lemma.} {\it If $(\PP , h)$ is soft, then
$(\PP ,h)$ has property (*).}
\smallskip
{\it Proof.} Suppose not and let 
$\{ p_n ;\; n \in \omega \}$ and $m \in \omega$
witness the contrary. For each $n \in \omega$ let $\{ q_j^n ;
\; j \in k_n \}$ be a weak finite cover with respect to
$\{ p_i ; \; i \in n \}$, $m$ according to (II). By assumption
none of these sets can be empty and we can assume 
that each $q_j^n$ has height $\leq m$. By the cover property (II) (ii)
they form an $\omega$-tree with finite levels with respect to "$\leq$". 
By K\"onig's Lemma this tree has an infinite branch.
By (I) there is a condition below this branch, contradicting the
fact that $\{ p_n ;\; n \in \omega \}$ is a maximal antichain. $\qed$
\bigskip
{\sanse 1.3} {\capit Theorem.} {\it Suppose $\PP$ is a ccc p.o.,
$h$ is a height function on $\PP$, and $(\PP , h)$ satisfies
property (*). Then any unbounded family of functions in $\omega^\omega
\cap V$ is still unbounded in $V[G]$, where $G$ is $\PP$-generic over $V$.}
\smallskip
{\it Proof.} Let $F$ be unbounded in $\omega^\omega \cap V$.
Suppose $\forces_{\PP} \breve f \in \omega^\omega$. For each $m \in \omega$
let $\{ p_n^m ; \; n \in \omega \}$ be a maximal antichain deciding the
value $\breve f(m)$. Choose $n_m$ according to (*) so that: whenever
$p$ is incompatible with $\{p_j^m ; j \in n_m \}$, then $h(p) > m$.
Define $f: \omega \to \omega$ by setting $f(m):=$ the maximum of the
values of $\breve f(m)$ decided by $\{p_j^m ; j \in n_m \}$. Let
$g \in F$ be a function which is not dominated by $f$.
We claim that $\forces_{\PP} $ "$\breve f$ does not dominate $g$".
\par For suppose there is a $p \in \PP$ and a $k \in \omega$ so that
\smallskip
\centerline{$p \forces \forall m \geq k \;\; \breve f (m) > g(m)$.}
\smallskip
\noindent Choose $m \geq k$ so that $h(p) \leq m$ and $g(m) \geq f(m)$.
Then $p$ must be compatible with $p_j^m$ for some $j \in n_m$.
But if $q$ is a common extension, then \smallskip
\centerline{$q \forces \breve f(m) > g(m) \geq f(m) \geq \breve f(m)$,}
\smallskip \noindent a contradiction. $\qed$
\bigskip
{\sanse 1.4} {\it Towards the proof of Theorem 1.} We think
of $\BB$ as consisting of sets $B \subseteq 2^\omega$
of positive measure so that for all $t \in 2^{<\omega}$,
if $[t] \cap B \neq \emptyset$ then $\mu ( [t] \cap
B ) > 0$; for $m \in \omega$ let $B \cap 2^m :=
\{ t \in 2^m ; \; [t] \cap B \neq \emptyset \}$. 
Then we define the 
following p.o. $(\PP,\leq)$:
\smallskip
\centerline{$(B,n) \in \PP \Longleftrightarrow B \in \BB \; \land
\; n \in \omega$;}
\smallskip
\centerline{$(B,n) \leq (C,m) \Longleftrightarrow B \subseteq C \; \land
\; n \geq m \; \land \; B \cap 2^m = C \cap 2^m$.}
\smallskip
\noindent 
It follows from the $ccc$-ness of any product of finitely many
copies [Je 2, part I, 5.7] of $\BB$ that $\PP$ is $ccc$, too. 
Clearly, $\PP$ generically adds a perfect set of random reals,
and we have to show that it does not add dominating reals. To this end,
we will introduce a height function on $\PP$. \par
In fact, let $\PP' \subseteq \PP$ be the set of conditions $(B,n) \in
\PP$ so that $\forall s \in 2^n \cap B \; \mu([s] \cap B) \geq
2^{-(lh(s) + 1)}$. $\PP'$ is dense in $\PP$ (by the Lebesgue density
Theorem [Ox, Theorem 3.20]). We define $h: \PP' \to \omega$ by $h((B,n)) = n$ and work
with $\PP'$ from now on.
\bigskip
{\sanse 1.5} {\capit Lemma.} {\it $(\PP' , h)$ is soft.}
\smallskip
{\it Remark.} By 1.2 and 1.3 the proof of this Lemma finishes the proof
of Theorem 1.
\smallskip
{\it Proof.} (I) is clear (for if $\{ (B_n , m) ; \;
n \in \omega \}$ is decreasing then $(\cap B_n ,m)$ is a lower
bound because we took our conditions from $\PP'$). 
For (II) we use:
\bigskip
{\sanse 1.6} {\capit Main Claim.} {\it Given $(B,n) , \; (C,m) \in
\PP'$ and $k \in \omega$ there are finitely many conditions $\{q_i ; \;
i \in j \}$ below $(C,m)$ so that \par
\item{(i)} each $q_i$ is incompatible with $(B,n)$; \par
\item{(ii)} if $q$ is incompatible with $(B,n)$, $h(q) \leq k$, and
$q \leq (C,m)$, then $\exists i \in j \; (q \leq q_i)$. \par}
\smallskip
{\it Proof.} Without loss $k \geq m, n$. Assume $n \geq m$. 
Let $\ell$ be such that $m \leq \ell \leq n$. We now describe which
conditions of height $\ell$ we put into our finite set. \par
(i)  For each $T \subseteq 2^\ell$ with $T \restrict m = C \cap 2^m$
and $T \subseteq C \cap 2^\ell$ and $T \neq B \cap 2^\ell$ let
$C_T \in \BB$ be such that $C_T \cap 2^\ell = T$ and $C_T \cap [t]
= C \cap [t]$ for each $t \in T$. If $(C_T, \ell) \in \PP'$, then
put $(C_T , \ell)$ into the set. \par
(ii) For each $T \subseteq 2^n$ with $T \restrict m = C \cap 2^m$
and $T \subseteq C \cap 2^n$ and $T \restrict \ell = B \cap 2^\ell$
and $T \not\supseteq B \cap 2^n$ let $C_T \in \BB$ 
be such that $C_T \cap
2^n =T$ and $ C_T \cap [t] = C \cap [t]$ for each $t \in T$.
If $(C_T, \ell) \in \PP'$, then put $(C_T , \ell)$ into the 
set. \par
(iii) For each $T \subseteq 2^n$ with $T \restrict m = C \cap 2^m$
and $T \subseteq C \cap 2^n$ and $T \restrict \ell = B \cap 2^\ell$
and $T \supseteq B \cap 2^n$ and for each $t \in B \cap 2^n$
let $C_{T,t} \in \BB$ be such that $C_{T,t} \cap 2^n = T$ and
$C_{T,t} \cap [s] = C \cap [s]$ for each $s \in T \setminus \{ t \}
$ and $C_{T,t} \cap [t] = (C \cap [t]) \setminus B$. If $(C_{T,t} , \ell)
\in \PP'$, then put $(C_{T,t}, \ell)$ into the set. \par
It is easy to see that any condition of height $\ell$ below $(C,m)$
which is incompatible with $(B,n)$ lies below one the conditions
defined in (i)  -- (iii) above. \par
Next suppose that $n \leq \ell \leq k$. Then we can again find a
finite set of conditions of height $\ell$ satisfying the requirements
of the main claim for conditions of height $\ell$ by an argument
similar to the one in (i) -- (iii) above. \par
This takes care of the case when $n \geq m$. So assume now $n \leq m$.
Then we get our set of conditions as in the preceding paragraph. $\qed$
\bigskip
{\sanse 1.7} {\it Proof of (II) of Lemma 1.5 from the main claim 1.6.} We make
induction using the main claim repeatedly. I.e. let $(B,n) = p_0$
and $(C,m) = (2^\omega , 0)$ and apply the main claim to them to
get $\{ q_i ; \; i \in j \}$. Then let $(B,n) = p_1$ and
$(C,m) = q_i$ $(i \in j)$ and apply the main claim $j$ times
to get a new family. Etc. $\qed$
\smallskip
This finishes the proof of Lemma 1.5 and of Theorem 1. $\qed$
\bigskip
{\sanse 1.8} {\capit Theorem.} {\it (i) Suppose $\langle \PP_\alpha,
(\breve Q_\alpha , \breve h_\alpha ) ; \; 
\alpha < \kappa \rangle$ is a finite support
iteration of arbitrary length $\kappa$ ($\kappa$ limit) such that
\smallskip
\centerline{$\forces_{\PP_\alpha} $ "$ \breve Q_\alpha$ is ccc, $\breve
h_\alpha$ is a height function on $\breve Q_\alpha$ and $(\breve Q_\alpha
, \breve h_\alpha)$ has property (*)".}
\smallskip
\noindent Then $\PP_\kappa = \lim_{\alpha < \kappa} \PP_\alpha$ does
not add dominating reals. \par
(ii) Suppose $\langle (\PP_\alpha , h_\alpha ) ; \; \alpha < \kappa \rangle$
is a sequence of soft ccc p.o.s of arbitrary length $\kappa$.
Then there is a height function $h$ on 
the finite support product $\PP$ of the $\PP_\alpha$ $(\alpha
< \kappa)$ so that $(\PP , h)$ is soft.} 
\smallskip
{\it Remark.} In particular, both the finite support iteration
and the finite support product of an arbitrary number of ccc
soft p.o.s does not add dominating reals (cf 1.2, 1.3).
\smallskip
{\it Proof.} (i) It suffices to show by induction on $\alpha$
that \smallskip
\centerline{$\forces_{\PP_\alpha} \omega^\omega \cap V$ is unbounded
in $\omega^\omega$.}
\smallskip
\noindent If $\alpha$ is a limit ordinal, this follows from
[JS 1, Theorem 2.2]. So suppose $\alpha$ is a successor.
Then the result follows from Theorem 
1.3 and the induction hypothesis.
\par
(ii) We make again induction on $\alpha$. Let $\QQ_\alpha$ be
the finite support product of the $\PP_\beta$ where $\beta
< \alpha$. We shall recursively construct height functions 
$g_\alpha : \QQ_\alpha \to \omega$ such that \par
\itemitem{(a)} for $\alpha < \beta$, $g_\alpha \subseteq g_\beta$;
\par
\itemitem{(b)} $g_\alpha (q) \geq \max_{\beta < \alpha} h_\beta (q 
\restrict \PP_\beta) $; \par 
\itemitem{(c)} $g_\alpha (q) \geq \vert supp (q) \vert$; \par
\itemitem{(d)} $(\QQ_\alpha , g_\alpha)$ is soft. \par
\noindent Clearly, a $g_\alpha$ satisfying (b) and (c) will satisfy
the decreasing chain property (I) as well. (We assume without
loss that $\forall \beta < \kappa$, $h_\beta (\11) = 0$.) \par
We first deal with the case when $\alpha$ is
a successor ordinal, $\alpha = \beta + 1$. Then $\QQ_\alpha =
\QQ_\beta \times \PP_\beta$. Let $m := \max \{ g_\beta (q),
h_\beta (p) \}$ and
define $g_\alpha :
\QQ_\alpha \to \omega$ by $$g_\alpha (q,p) := 
\cases{m+1 & if $\vert supp (q) \vert = m$ and $p \neq \11$ \cr
m & otherwise}$$
for $(q,p) \in \QQ_\beta \times \PP_\beta$.
$g_\alpha$ is a height function on $\QQ_\alpha$ which is
easily seen to satisfy  (a) --- (c) above.
\par
To show that $(\QQ_\alpha, g_\alpha)$ satisfies the weak finite
cover property (II), let $\{ (q_i , p_i) ; \; i \in n \}$
be a finite subset of $\QQ_\alpha$ and let $m \in \omega$.
For each $A \subseteq n$ let $\{ q_j^A ; \; j \in k_A \}$
be a weak finite cover with respect to $\{ q_i ; \; i \in A \}$
and $m$ in $\QQ_\beta$ (i.e. (i) $\forall i \in A, \; j \in k_A$,
$q_j^A$ is incompatible with $q_i$; and (ii) whenever $q$ is
incompatible with all $q_i$ ($i \in A$) and $h(q) 
\leq m$ then there exists $j \in k_A$ so that $q \leq q_j^A$),
and let $\{ p_j^A ; \; j \in \ell_A \}$ be a weak finite
cover with respect to $\{ p_i ; \; i \in n \setminus A \}$ and $m$.
We claim that the family $F :=\{ (q_i^A , p_j^A ) ; \; A \subseteq n \;
\land \; i \in k_A \; \land\; j \in \ell_A \}$ is a weak
finite cover with respect to $\{ (q_i , p_i) ; \; i \in n \}$ and $m$.
\par
Clearly, $F$ satisfies (i) of the definition of the weak finite
cover property (II). Furthermore, if $(q,p)$ is incompatible with
all $(q_i , p_i)$ ($i \in n$) there exists $A \subseteq n$ such
that $q$ is incompatible with all $q_i$ for $ i \in A$ and $p$ is
incompatible with all $p_i$ for $i \in n \setminus A$. So if
$g_\alpha (q,p) \leq m$ (in particular, $g_\beta (q) \leq
m$ and $h_\beta (p) \leq m$) then we can find $j \in k_A$ and $j'
\in \ell_A$ such that $q \leq q_j^A$ and $p \leq p_{j'}^A$;
i.e. $(q,p) \leq (q_j^A , p_{j'}^A)$. This shows (ii) in the
definition of the weak finite cover property (II).
\par
Now suppose $\alpha$ is a limit ordinal. Then let $g_\alpha :=
\bigcup_{\beta < \alpha} g_\beta$. $g_\alpha$ clearly satisfies
(a) --- (c), and the weak finite cover property (II) 
for $(\QQ_\alpha, g_\alpha)$ follows from the weak finite cover
properties of the $(\QQ_\beta , g_\beta)$ for $\beta < \alpha$
(because (II) talks only about finitely many conditions).
$\qed$
\bigskip
{\sanse 1.9} 
We note here that the notions discussed so far are closely
tied up with some cardinal invariants of the continuum.
Namely, we let ${\cal N}$ denote the ideal of null sets and
\smallskip
\itemitem{$add({\cal N}) :=$} the least $\kappa$ such that 
$\exists {\cal F} \in [{\cal N}]^\kappa \; (\bigcup
{\cal F} \not\in {\cal N})$;
\par
\itemitem{$wcov ({\cal N}) :=$} the least $\kappa$ such that
$\exists {\cal F} \in [{\cal N}]^\kappa \; (2^\omega \setminus\bigcup
{\cal F}$ does not contain a perfect set);
\par
\itemitem{$cov ({\cal N}) :=$} the least $\kappa$ such that
$\exists {\cal F} \in [{\cal N}]^\kappa \; (\bigcup
{\cal F} = 2^\omega)$;
\par
\itemitem{$unif({\cal N}) :=$} the least $\kappa$ such that
$[2^\omega]^\kappa \setminus {\cal N} \neq \emptyset$;
\par
\itemitem{$wunif ({\cal N}) :=$} the least $\kappa$ such that 
there is a family ${\cal F} \in [[2^{< \omega}]^\omega]^\kappa$
of perfect sets with $\forall N \in {\cal N} \; \exists
T \in {\cal F} \; (N \cap T = \emptyset)$;
\par
\itemitem{$cof ({\cal N}) :=$} the least $\kappa$ such that
$\exists {\cal F} \in [{\cal N}]^\kappa \; \forall A \in {\cal N}
\; \exists B \in {\cal F} \; (A \subseteq B)$;
\par
\itemitem{$b :=$} the least $\kappa$ such that $\exists {\cal F}
\in [\omega^\omega]^\kappa \; \forall f \in \omega^\omega
\; \exists g \in {\cal F} \; \exists^\infty n \;
(g(n) > f(n))$;
\par
\itemitem{$d :=$} the least $\kappa$ such that $\exists {\cal F}
\in [\omega^\omega]^\kappa \; \forall f \in \omega^\omega \;
\exists g \in {\cal F} \; \forall^\infty n \; (g(n) > f(n))$.
\smallskip
\noindent Then we can arrange these cardinals in the following diagram.
\bigskip
\centerline{put diagram 1 here}
\bigskip
\noindent Here the invariants get larger as one moves up in the
diagram. $b \geq add({\cal N})$ (and dually $d \leq cof({\cal N})$)
is due to Miller [Mi]. The dotted line says that
$wcov ({\cal N}) \geq min \{ cov({\cal N}) , b \}$ (and dually,
$wunif ({\cal N}) \leq max \{ unif ({\cal N}) , d \}$).
This can be seen from the Bartoszy\'nski--Judah result (*) 
in the Introduction as follows. Suppose $\lambda := wcov 
({\cal N}) < min \{ cov({\cal N}) , b \}$.
Let $M$ be a model of enough $ZFC$ of size $\lambda$ containing
a weak covering family. As $\lambda < b$ there is a real 
$f \in \omega^\omega$ dominating all reals in $M$. Let $N$ be a model 
of enough $ZFC$ of size $\lambda$ containing $M$ and $f$. As $\lambda
< cov ({\cal N})$, there is a real $r \in 2^\omega$ random over
$N$. By (*) this implies that there is a perfect set of random
reals over $M$, a contradiction. --- Iterating the p.o. from Theorem
1 we get:
\smallskip
{\capit Theorem 1'.} {\it For any regular cardinal $\kappa$, it is
consistent that $wcov({\cal N}) = \kappa$ while $b = \omega_1$;
Dually, it is consistent that $wunif ({\cal N}) = \omega_1$ while
$d = \kappa$.}
\smallskip
{\it Proof.} (a) Assume $CH$. We make a finite support
iteration of length $\kappa$ of the p.o. $\PP$ described in 1.4.
In the generic extension we have $wcov ({\cal N}) = \kappa$ because
we added $\kappa$ many perfect sets of random reals; and $b = \omega_1$
by 1.5, 1.2 and 1.8 (i). 
\par
(b) Assume $MA + 2^\omega = \kappa$; and make a finite
support iteration of length $\omega_1$ of $\PP$. Again
standard arguments show that $wunif({\cal N}) = \omega_1$
and $d = \kappa$ in the generic extension. $\qed$
\smallskip
The most interesting open question concerning the
relationship between these cardinals is connected with
Question 3 in the Introduction.
\smallskip
{\dunh Question 3'.} {\it Is it consistent that $wcov({\cal N})
> d$? Dually, is it consistent that $wunif({\cal N}) < b$?}
\bigskip
{\sanse 1.10} {\it Proof of Theorem 2.} By 1.2, 1.3 and 1.8 (ii)
it suffices to show that there is a height function $h : \BB \to
\omega$ so that $(\BB , h)$ is soft. But this is easy: for $B \in
\BB$ let $h(B) := \min \{ n \in \omega ; \; \mu (B) \geq {1 \over n} \}$.
$\qed$
\smallskip
We note that this height function $h$ and also
the height function it induces on $\BB \times \BB$ by 1.8 (ii)
have a {\it strong} finite cover property: (ii) in (II)
can be replaced by: whenever $q$ is incompatible with all
$p_i$ then there exists $j \in k$ so that $q \leq q_j$.

\vfill\eject

\noindent{\dunhg $\S$ 2 Not adding perfect sets of random reals}
\Smallskip
{\sanse 2.1} This whole section is devoted to the proof of Theorem 3.
Lemmata 2.2 and 2.3 below which we single out from the principal
argument bear
the imprint of the proof of Cicho\'n's Theorem in [BJ 1, 2.1 -- 2.4].
The main new idea comes in in 2.4. The rest (2.5 -- 2.9) is mostly
technical.
\par
Given $k', k \in \omega$, $k' < k$, we let 
$\epsilon_{k,k'} := 2^{1-k} \cdot (1 + {k \choose 1} + ... +{k \choose
k'-1})$. Clearly, given any $k' \in \omega$, we can find $k > k'$
so that $\epsilon_{k,k'}$ is arbitrarily small. 
\bigskip
{\sanse 2.2} {\capit Lemma.} {\it Given $n, k, k' \in \omega$ ($k' \leq k$ and
$k \leq 2^n$) and $Z \subseteq 2^n$ and real numbers $a_T$ for
each $T \subseteq 2^n$ with $\vert T \cap Z \vert \geq k$ there
exists $Z' \subseteq Z$ of size $\leq {\vert Z \vert \over 2}$
such that $\sum_{\vert T \cap Z' \vert \geq k'} a_T \geq
\sum_T a_T \cdot (1 - \epsilon_{k,k'})$.}
\smallskip
{\it Proof.} Let $a:=\sum_T a_T$. For any $\lambda$ close to $1$
($\lambda < 1$) we can choose $\ell \in \omega$ and $\{ T_i ; \;
i \in \ell \}$ so that 
$$a_T \cdot \lambda < {\vert \{ i ; \; T_i = T \} \vert \over \ell}
\cdot a < a_T \cdot \lambda^{-1}$$
for all $T$. For $i \in \ell$ let ${\cal Z}_i := \{ Z' \subseteq Z ;
\; \vert Z' \cap T_i \vert \geq k'$ and $\vert (Z \setminus
Z') \cap T_i \vert \geq k' \}$.
Then $\vert {\cal Z}_i \vert \cdot 2^{- \vert Z \vert} \geq
1 - \epsilon_{k,k'}$ for all $i \in \ell$. We claim that
there is an $X \subseteq \ell$ of size $\geq \ell \cdot(1-
\epsilon_{k,k'})$ so that $\bigcap_{i \in X} {\cal Z}_i \neq
\emptyset$. \par
For suppose not. Then for each $Z' \subseteq Z$, $\vert \{ i \in
\ell
; \; Z' \in {\cal Z}_i \} \vert < \ell \cdot (1 - \epsilon_{k,k'})$.
Hence $2^{\vert Z \vert} \cdot 
\ell \cdot (1 - \epsilon_{k,k'}) \leq \sum_{i \in \ell}
\vert {\cal Z}_i \vert = \sum_{Z' \subseteq Z} \vert \{ i 
\in \ell ; \; Z' \in {\cal Z}_i \} \vert < 2^{\vert Z \vert} \cdot
\ell \cdot (1 - \epsilon_{k,k'})$, a contradiction. \par
Now choose $Z' \in \bigcap_{i \in X} {\cal Z}_i$. Then either
$\vert Z' \vert \leq {\vert Z \vert \over 2}$ or $\vert Z \setminus Z'
\vert
\leq {\vert Z \vert \over 2}$. Assume without loss that $\vert
Z' \vert \leq {\vert Z \vert \over 2}$. Furthermore
$$\sum_{\vert T \cap Z' \vert \geq k'} a_T \cdot \lambda^{-1} >
\sum_{\vert T \cap Z' \vert \geq k'} {\vert \{ i ; \; T_i = T \}
\vert \over \ell} \cdot a \geq {\vert X \vert \over \ell} \cdot a
\geq a \cdot (1 - \epsilon_{k,k'}).$$
Because there are only finitely many possibilities for the sum on
the lefthand side, we can choose $\lambda$ so small that for the
$Z'$ chosen according to this $\lambda$ we have
$$\sum_{\vert T \cap Z' \vert \geq k'} a_T \geq a \cdot (1 -
\epsilon_{k,k'}).$$
This finishes the proof of the Lemma. $\qed$
\bigskip
{\sanse 2.3} {\capit Lemma.} {\it Given  a real $\epsilon > 0$ and $m \in
\omega$ the following is true for large enough $k, n \in \omega$:
given real numbers $a_T$ for each $T \subseteq 2^n$ with
$\vert T\vert \geq k$ there exists a $Z \subseteq 2^n$ of size
$\leq 2^{n - m}$ such that $\sum_{T \cap Z \neq \emptyset}
a_T \geq \sum_T a_T \cdot (1 - \epsilon)$.}
\smallskip
{\it Remark.} We say that a statement is true for {\it large
enough $n$} iff $\exists k \in \omega$ so that $\forall n \geq
k$ the statement is true.
\smallskip
{\it Proof.} Construct recursively a sequence $\langle k_i ; \;
i \leq m \rangle$ of natural numbers so that $\prod_{i \in m}
(1 - \epsilon_{k_{i+1},k_i}) \geq (1 - \epsilon)$ where $k_0 = 1$. Let
$k \geq k_m$ and $n$ so large that $k \leq 2^n$.
Now apply Lemma 2.2 $m$ times to get $Z$. $\qed$
\bigskip
{\sanse 2.4} {\it Diagonal chains.} It turns out that 
a detailed investigation of antichains in $\BB \times
\BB$ is necessary for the proof of Theorem 3. 
We say $(p,q) \in \BB \times \BB$ is {\it quadratic} iff
$\mu(p) = \mu(q)$. Clearly the quadratic conditions are dense
in $\BB \times \BB$ so that it suffices (essentially) to consider
them. --- More generally, given $(p,q) \in \BB \times \BB$,
$(p',q') $ is {\it quadratic in $(p,q)$} iff
$p' \leq p, \; q' \leq q$ and ${\mu(p') \over \mu(p)} =
{\mu(q') \over \mu(q)}$. \par
$\{ (p_n, q_n ) ; \; n \in \omega \}$ is said to be a
{\it first order diagonal chain} in $\BB \times \BB$ iff \par
\itemitem{(1)} each $(p_n, q_n)$ is quadratic; \par
\itemitem{(2)} both $\{ p_n ; \; n \in \omega \}$ and $\{ q_n ; \;
n \in \omega \}$ are maximal antichains in $\BB$. \par
\noindent More generally we say that $C=\{(p_n^{\sigma\tau}, q_n^{\sigma
\tau} ) ; \; n \in \omega, \; \sigma, \tau \in \omega^{<m}, \; 
lh(\sigma)= lh(\tau), \; \forall i \in lh(\sigma) \;
(\sigma (i) \neq \tau (i)) \}$ is an {\it $m$th order diagonal
chain} in $\BB \times \BB$ (for $m \geq 2$)
iff 
$\{ ( p_n^{\sigma\tau}, q_n^{\sigma\tau}) \in C ; \;
lh(\sigma) < m-1 \}$ is an $(m-1)$th order diagonal chain and
for each $\sigma, \; \tau$ of length $m-1$ with $\tau (m-2)
\neq \sigma(m-2)$, \par  
\itemitem{(1)} $\forall n \in \omega$, $(p_n^{\sigma\tau}, q_n^{\sigma\tau})$
is quadratic in $(p_{\sigma(m-2)}^{\sigma \restrict (m-2) \tau
\restrict (m-2)} , q_{\tau(m-2)}^{\sigma \restrict (m-2) \tau
\restrict (m-2)})$; \par
\itemitem{(2)} $\{ p_n^{\sigma\tau} ; \;
n \in \omega 
\}$ is a maximal antichain of conditions below
$p_{\sigma (m-2)}^{\sigma \restrict (m-2) \tau \restrict (m-2)}$
in $\BB$ and $\{ q_n^{\sigma\tau} ; \; n \in \omega \}$ is a 
maximal antichain of conditions below $q_{\tau(m-2)}^{\sigma
\restrict (m-2) \tau \restrict (m-2)}$ in $\BB$. \par
\bigskip
\centerline{put diagram 2 here}
\bigskip
Clearly, below any antichain $A$ in $\BB \times \BB$ there is an
$m$th order diagonal chain $C$ for each $m$ (in the sense that
any condition in $C$ is below some condition in $A$).
\bigskip
{\sanse 2.5} {\it Towards the proof of Theorem 3.}
Let $\breve T$ be a $\BB \times \BB$-name so that
\smallskip
\centerline{$\forces_{\BB \times \BB} $ "$\breve T$ is perfect".}
\smallskip
\noindent We want to construct a null set $N$ in the ground model so that
\smallskip
\centerline{$\forces_{\BB \times \BB} $ "$N \cap \breve T \neq \emptyset$".}
\smallskip
\noindent More explicitly, using Lemma 2.3, we shall construct sequences
$\langle n_m ; \; m \in \omega \rangle$ and $\langle Z_m ; \;
m \in \omega \rangle$ so that $Z_m \subseteq 2^{n_m}$, $\vert Z_m
\vert \leq 2^{n_m - m}$ and \smallskip
\centerline{$\forces_{\BB \times \BB} \exists x \in \breve T \;
\exists^\infty m \; ( x \restrict n_m \in Z_m )$.}
\smallskip
\noindent This will imply the required result for $N := \{
x \in 2^\omega ; \; \exists^\infty m \; (x \restrict n_m \in Z_m) \}
$ is a null set (see below in 2.9). --- We set $n_0 := 0$ and
$Z_0 := \emptyset$. Now let $m > 0$ and assume that $n_{m-1}$
and $Z_{m-1}$ have been defined. 
We shall describe the construction
of $n_m$ and $Z_m$.
\smallskip
Let $\delta > 0$ be very small; let $\langle z_j ; \; j \in m \rangle$,
$\langle y_j ; \; j \in m - 1 \rangle$ be sequences of natural numbers so
that $z_0 > 1$, $y_j = 4 \cdot z_j$ and $y_j \over z_{j + 1}$ is very
small; 
let $\epsilon > 0$ such that $\delta^{-2} \cdot \epsilon \cdot
m \cdot z_{m - 1}$ is very small (in fact, we want that $\zeta_m
= \zeta = 2 \cdot m \cdot ( \epsilon + \delta^{-2} \cdot \epsilon \cdot
m \cdot z_{m-1} + \delta ) + \sum_{j = 0}^{m - 2} {y_j \over 
z_{j + 1}}$ is -- say -- smaller than $1 \over 4^m$ -- cf 2.6);
let $v \in \omega$ be such that $z_{m-1} \leq v \cdot (1 -
\epsilon)^2$.
Choose $k$ according to Lemma 2.3 for $\epsilon$ and $m + 2 \cdot 
n_{m-1}$.
Let $\{ (p_i^{\sigma\tau}, q_i^{\sigma\tau}); \; i \in \omega,
\; \sigma, \tau \in \omega^{<m}, \; lh(\sigma) = lh(\tau) , \;
\forall j \in lh(\sigma) \;(\sigma(j) \neq \tau(j)) \}$ be an
$m$th order diagonal chain in $\BB \times \BB$ deciding $\breve T$
up to some level $n_i^{\sigma\tau}$ such that for each $s \in
2^{\leq n_{m-1}}$,
\smallskip
\centerline{either $(p_i^{\sigma\tau}, q_i^{\sigma\tau}) \forces 
s \in \breve T \; \land \; \vert \breve T_s
\restrict n_i^{\sigma\tau} \vert \geq k$}
\smallskip
\centerline{or
$(p_i^{\sigma\tau}, q_i^{\sigma\tau}) \forces s \not\in \breve T$.}
\smallskip
\noindent We can construct this diagonal chain in such a way that
for $\sigma, \tau$ of length $\ell < m$ (with $\forall j \in \ell
\; ( \sigma (j) \neq \tau (j))$),
$$\sum_i \mu(p_i^{\sigma\tau} , q_i^{\sigma\tau} ) < {1 \over v}
\cdot \mu(p_{\sigma (\ell - 1)}^{\sigma \restrict (\ell - 1) \tau
\restrict (\ell - 1)}, q_{\tau (\ell - 1)}^{\sigma \restrict
(\ell - 1) \tau \restrict (\ell - 1)}) \eqno (5.1)$$
(where we make the convention that for $\ell = 0$, $\sigma =
\tau = \langle\rangle$,
$p_{\sigma (\ell - 1)}^{\sigma \restrict (\ell - 1) \tau
\restrict (\ell - 1)} = q_{\tau (\ell - 1)}^{\sigma \restrict
(\ell - 1) \tau \restrict (\ell - 1)} = \11$,
the maximal element of $\BB$).
\par
We now define recursively which pairs of sequences
$\sigma\tau$ ($\sigma , \tau \in \omega^{< m}$, $lh(\sigma)
= lh(\tau)$, $\forall j \in lh(\sigma) \; (\sigma (j) \neq
\tau (j))$) are {\it relevant}. For relevant pairs we also
define $a^{\sigma\tau} \in \RR$ and $j^{\sigma\tau} \in
\omega$. $\langle\rangle\langle\rangle$ is relevant. 
Choose $j^{\langle\rangle\langle\rangle} \in \omega$ such that
$a^{\langle\rangle\langle\rangle} := \sum_{i \in j^{\langle\rangle\langle
\rangle}} \mu(p_i^{\langle\rangle\langle\rangle}) \geq 1 - \epsilon$.
Suppose $a^{\sigma\tau}$ and $j^{\sigma\tau}$ are defined for relevant
pairs $\sigma\tau$ of length $\ell$ ($0 \leq \ell < m-1$). Then
$\sigma \hat{\;}\langle i \rangle \tau \hat{\;} \langle j \rangle$
is relevant iff $i,j \in j^{\sigma\tau}$ and $i \neq j$. 
Furthermore, for each such $i,j$, choose $j^{\sigma\hat{\;}\langle i \rangle
\tau \hat{\;} \langle j \rangle} \in \omega$ such that
$$a^{\sigma\hat{\;} \langle i \rangle \tau \hat{\;} \langle j \rangle}
:= \sum_{i' \in j^{\sigma \hat{\;} \langle i \rangle\tau \hat{\;} \langle j
\rangle}} \mu(q_{i'}^{\sigma\hat{\;}\langle i \rangle \tau\hat{\;}
\langle j \rangle}) \cdot \mu (p_i^{\sigma\tau}) \geq \mu (p_i^{\sigma
\tau}) \cdot \mu(q_j^{\sigma\tau}) \cdot ( 1 - \epsilon) . \eqno
(5.2) $$
Now let $$n_m := n:= \max_{\sigma \tau \;{\rm relevant}, i \in j^{\sigma
\tau}} n_i^{\sigma\tau}.$$
\par
Fix $s \in 2^{\leq n_{m-1}}$. For $T \subseteq 2^n$ 
with $\vert T \vert \geq k$ and $s \subseteq stem(T)$ and for relevant tuples
$\sigma
\tau$ let $$a_T^{\sigma\tau} := \sum_{i \in j^{\sigma\tau}} {\mu
(\parallel s \in \breve T \; \land \;
\breve T_s \restrict n = T \parallel \cap (p_i^{\sigma\tau},
q_i^{\sigma\tau})) \over \mu(p_i^{\sigma\tau})} \cdot \mu(p_{\sigma
(\ell - 1)}^{\sigma \restrict (\ell - 1 )\tau \restrict (\ell - 1)}).
\eqno (5.3) $$
And let
$$a_s^{\sigma\tau} := \sum_{i \in j^{\sigma\tau}} { \mu
(\parallel s \not\in \breve T \parallel \cap (p_i^{\sigma\tau},
q_i^{\sigma\tau})) \over \mu (p_i^{\sigma\tau})} \cdot \mu
(p_{\sigma (\ell - 1)}^{\sigma \restrict ( \ell - 1 ) \tau
\restrict (\ell - 1)} ).$$
Then $\sum_T a_T^{\sigma\tau} + a_s^{\sigma\tau} = a^{\sigma\tau}$ (this uses
the finite additivity of the measure on $r.o.(\BB \times \BB)$ --
see Introduction). Let $a_T^j :=
\sum_{lh(\sigma) = lh (\tau) = j} a_T^{\sigma\tau}$,
$a_s^j := \sum_{lh(\sigma) = lh(\tau) = j} a_s^{\sigma\tau}$
and $a^j := \sum_T a_T^j  + a_s^j = \sum_{lh(\sigma) = lh(\tau) = j} a^{
\sigma\tau}$. Let $a_T := \sum_{j \in m} {a_T^j \over a^j -
a_s^j}$.
Apply Lemma 2.3 to get $Z_m^s := Z^s \subseteq 2^n$ of 
size $\leq 2^{n - m - 2 \cdot n_{m-1}}$
such that $$\sum_{T \cap Z^s \neq \emptyset} a_T \geq \sum_T a_T \cdot
(1 - \epsilon) = m \cdot ( 1 - \epsilon) . \eqno (5.4) $$
\par
Finally, set $Z_m := Z:= \bigcup_{s \in 2^{\leq n_{m-1}}} Z^s$.
This completes the construction of $n_m = n$ and $Z_m = Z$.
\bigskip
{\sanse 2.6} {\capit Main Claim.} {\it Let $s \in
2^{\leq n_{m-1}}$. Suppose $(p,q)$ is a quadratic condition
such that $(p,q) \forces $ "$ s \in \breve T \land Z_m^s \cap 
\breve T_s = \emptyset$". Then
$\mu (p,q) < {1 \over 4^m} + \zeta_m$ (where $\zeta_m = \zeta =
2 \cdot m \cdot ( \epsilon + \delta^{-2} \cdot m \cdot z_{m-1}
+ \delta) + \sum_{j = 0}^{m-2} {y_j \over z_{j+1}}$ as in 2.5).}
\smallskip
{\it Proof.} The proof of the main claim will take some time (up to 2.8);
to make our argument (which is essentially one big estimation) go
through smoothly we need to make some conventions and 
introduce a few more notions. \par 
If $\sigma$ is a sequence of length $\ell \geq 1$, $\hat\sigma
= \sigma \restrict (\ell - 1)$ will be the sequence with the last
value deleted. For $\ell = 0$, $p_{\langle\rangle ( \ell - 1)}^{
\hat{\langle\rangle}\hat{\langle\rangle}} = q_{\langle\rangle
(\ell - 1)}^{\hat{\langle\rangle}\hat{\langle\rangle}} = \11$,
the maximal element of the Boolean algebra $\BB$. $\sum_{ij}$
will always stand for $\sum_{i,j \in j^{\sigma\tau}, i \neq j}$,
where $\sigma\tau$ is clear from the context; similarly, $\sum_{
\sigma\tau}$ means that the sum runs over all relevant $\sigma\tau$
of some fixed length $\ell$ (where $\ell$ is again clear
from the context). ---
For a relevant pair $\sigma\tau$ we let
$$\eqalignno{
A^{\sigma\tau} & := \{ i \in j^{\sigma\tau} ; \; \mu(p \cap 
p_i^{\sigma\tau} ) < \delta \cdot \mu(p_i^{\sigma\tau}) \} & (6.1)
\cr
B^{\sigma\tau} & := \{ i \in j^{\sigma\tau} ; \; \mu(q
\cap q_i^{\sigma\tau} ) < \delta \cdot \mu(q_i^{\sigma\tau}) \} \setminus
A^{\sigma\tau} & (6.2) \cr
C^{\sigma\tau} & := j^{\sigma\tau} \setminus ( A^{\sigma\tau}
\cup B^{\sigma\tau} ) & (6.3) \cr}$$
Let $lh(\sigma) = lh(\tau) = m -1$. We say the relevant
pair $\sigma\tau$ is {\it nice} iff 
$$\sum_{i \in C^{\sigma\tau}} \mu(q_i^{\sigma\tau}) \cdot \mu(
p_{\sigma(m-2)}^{\hat\sigma \hat\tau  })
\leq \delta^{-2} \cdot a^{\sigma\tau} \cdot \epsilon \cdot m \cdot 
z_{m-1}. \eqno (6.4) $$
More generally, if $lh(\sigma) = lh(\tau) = \ell$ (where
$0 \leq \ell < m-1$), we say the pair $\sigma\tau$ is {\it nice} iff
$${\rm (I)} 
\sum_{i \in C^{\sigma\tau}} \mu(q_i^{\sigma\tau}) \cdot \mu(
p_{\sigma(\ell - 1)}^{\hat\sigma  \hat\tau 
})
\leq \delta^{-2} \cdot a^{\sigma\tau} \cdot \epsilon \cdot m \cdot 
z_{m-1}; \eqno (6.5) $$
$${\rm (II)} \sum_{ij  , 
\sigma \hat{\;} \langle i \rangle \tau \hat{\;}
\langle j \rangle {\rm nice}} a^{\sigma \hat{\;} \langle i \rangle
\tau \hat{\;} \langle j \rangle} \geq (1 - {y_\ell \over z_{\ell + 1} })
\cdot
\sum_{ij 
} a^{\sigma \hat{\;} \langle i \rangle \tau \hat{\;}
\langle j \rangle}. \eqno (6.6) $$
(Note that this is a definition by backwards recursion on $\ell$.)
\bigskip
{\sanse 2.7} {\capit Claim.} {\it for any $\ell$ $(0 \leq \ell \leq m-1)$, 
$\sum_{ \sigma\tau  \;
{\rm nice}} a^{\sigma\tau} \geq (1 -{1 \over z_\ell}   )
\cdot a^\ell$.} \hfill{(7.1)}
\smallskip
{\it Proof.} 
We first show that for any $\ell$ we have
$$ 
\delta^{-2} \cdot a^\ell \cdot \epsilon \cdot m \geq \sum_{\sigma\tau} 
\sum_{i \in C^{\sigma\tau}} \mu(q_i^{\sigma\tau}) \cdot\mu(p_{\sigma
(\ell
- 1)}^{\hat\sigma  \hat\tau }). \eqno (7.2)$$
\par
By construction (5.4), $\sum_{T \cap Z^s \neq \emptyset} a_T \geq
m\cdot (1 - \epsilon)$; i.e. $\sum_{T \cap Z^s \neq
\emptyset} \sum_{j \in m} {a_T^j \over a^j - a_s^j} \geq m \cdot (1 -
\epsilon)$. Hence $\sum_{T \cap Z^s \neq \emptyset} a_T^\ell
\geq (a^\ell - a_s^\ell) (1 - \epsilon \cdot m)$. As $(p,q) \forces
"s \in \breve T \land \breve T_s \restrict n 
\cap Z^s = \emptyset"$ we get (using (5.3) and also the definition
of $C^{\sigma\tau}$ ((6.1) --- (6.3)))
$$\eqalign{a^\ell \cdot \epsilon \cdot m \geq 
(a^\ell - a^\ell_s) \cdot \epsilon \cdot m \geq & \sum_{\sigma\tau}
\sum_{i \in C^{\sigma\tau}} \mu(p_{\sigma(\ell - 1)}^{
\hat\sigma \hat\tau }) \cdot
{\mu(p \cap p_i^{\sigma\tau} , q \cap q_i^{\sigma\tau})
\over \mu (p_i^{\sigma\tau})} \cr \geq & \delta^2 \cdot\sum_{\sigma\tau}
\sum_{i \in C^{\sigma\tau}} \mu(q_i^{\sigma\tau})
\cdot\mu(p_{\sigma(\ell - 1)}^{\hat\sigma  \hat\tau 
}).\cr}$$
This shows that formula (7.2)
holds. \par
Next, we prove the claim by backwards induction. So assume
$\ell = m-1$. In that case, it follows immediately from formula
(7.2) and the definition of niceness (6.4) 
that $\sum_{\sigma\tau \;{\rm nice}} a^{\sigma\tau} \geq (1 -
{1 \over z_{m-1}}) \cdot a^{m - 1}$. \par
So let $\ell < m-1$ and assume the claim has been proved for $\ell
+ 1$. We let $\Sigma (I) : = \{ \sigma\tau;
\; lh(\sigma) = lh(\tau) = \ell$ and $\sigma\tau$ satisfies (I)
of the definition of niceness $\}$ and $\Sigma (II) :=
\{ \sigma\tau ; \; lh(\sigma) = lh(\tau) = \ell$ and
$\sigma\tau$ satisfies (II) of the definition of
niceness $\}$. By the argument of the preceding paragraph we know that
$\sum_{\sigma \tau \in \Sigma (I)} a^{\sigma\tau} \geq
(1 - {1 \over z_{m-1}}) \cdot a^\ell$. We claim that $$\sum_{\sigma\tau
\in \Sigma (II) } \sum_{ij 
} a^{\sigma \hat{\;} \langle i
\rangle \tau \hat{\;} \langle j \rangle} \geq ( 1 - {1 \over y_\ell})
\cdot a^{\ell +1}. \eqno (7.3) $$ \par
For suppose not. Then $\sum_{\sigma \tau \not\in \Sigma (II)
} \sum_{ij } a^{\sigma 
\hat{\;} \langle i \rangle \tau \hat{\;}
\langle j \rangle} > {1 \over y_\ell} \cdot a^{\ell +1}$. But if $\sigma\tau$
does not satisfy (II) then $$\sum_{ij 
, \sigma\hat{\;}\langle i \rangle
\tau\hat{\;}\langle j \rangle {\rm not \; nice}} a^{\sigma \hat{\;} \langle
i \rangle \tau \hat{\;} \langle j \rangle} > {y_\ell \over z_{\ell + 1}}
\cdot\sum_{ij } 
a^{\sigma \hat{\;} \langle i \rangle \tau \hat{\;} \langle j
\rangle}.$$ Hence $$\sum_{\sigma\tau \not\in \Sigma (II)}
\;\;\;\sum_{ij , \sigma \hat{\;} 
\langle i \rangle \tau \hat{\;} \langle j \rangle
{\rm not \; nice}} a^{\sigma \hat{\;} \langle i \rangle \tau \hat{\;}
\langle j \rangle} > {1 \over z_{\ell + 1}} \cdot a^{\ell + 1},$$ contradicting
the induction hypothesis.\par
As $(1 - \epsilon) \cdot a^\ell \leq {1 \over 1 - \epsilon} \cdot a^{\ell + 1} +
\sum_{\sigma\tau} 
\sum_{i \in j^{\sigma\tau}} \mu(p_i^{\sigma\tau}, q_i^{\sigma\tau})
< {1 \over 1 - \epsilon} \cdot a^{\ell + 1} + {1 \over v}
\cdot {1 \over 1 - \epsilon} \cdot a^\ell$ ((5.1) and (5.2)),
we have ${1 \over (1 - \epsilon)^2} \cdot 
a^{\ell + 1} \geq  ( 1 - {1 \over
v} \cdot {1 \over (1 - \epsilon)^2}) 
\cdot a^\ell$. Hence by (7.3) $$\sum_{\sigma \tau \in \Sigma (II)
} a^{\sigma\tau} \geq \sum_{\sigma\tau \in \Sigma (II)
} \sum_{ij } 
a^{\sigma \hat{\;} \langle i \rangle \tau \hat{\;} \langle j \rangle}
\geq (1 - {1 \over y_\ell}) \cdot a^{\ell + 1} 
\geq (1 - {1 \over y_\ell} - {1 \over v}
\cdot {1 \over (1 - \epsilon)^2} - 2 \cdot \epsilon) \cdot a^\ell.$$
Putting everything together we get that $$\sum_{\sigma\tau \;{\rm
nice} } a^{\sigma\tau} 
\geq (1 - {1 \over z_{m-1}} - {1 \over y_\ell} 
- {1 \over v} \cdot {1 \over (1 - \epsilon)^2} - 2 \cdot \epsilon ) 
\cdot a^\ell \geq (1 - {1 \over z_\ell}) \cdot a^\ell.
\;\;\;\qed$$
\bigskip
This shows in particular that the pair
$\langle\rangle\langle\rangle$
is nice.
\bigskip
{\sanse 2.8} {\capit Claim.} {\it If $\sigma\tau$ is nice of length $\ell$
$(0 \leq \ell \leq m-1)$ then $$\mu(p \cap p_{\sigma(\ell 
- 1)}^{\hat\sigma  \hat\tau },
q \cap q_{\tau(\ell - 1)}^{\hat\sigma \hat\tau
}) < ({1 \over 4^{m - \ell}} + \zeta_\ell) \cdot
\mu(p_{\sigma(\ell - 1)}^{\hat\sigma \hat\tau 
} , q_{\tau(\ell - 1)}^{\hat\sigma \hat\tau
}), \eqno (8.1)$$ where $\zeta_\ell := 2 \cdot (m-\ell) \cdot
(\epsilon + \delta^{-2} \cdot \epsilon \cdot m \cdot z_{m-1} +
\delta) + \sum_{j=\ell}^{m-2} {y_j \over z_{j+1}}$ $($in
particular $\zeta_0 = \zeta)$.}
\smallskip
{\it Proof.} 
We know that for arbitrary nice $\sigma\tau$, $$\delta^{-2} \cdot
a^{\sigma\tau} \cdot \epsilon \cdot m \cdot z_{m-1} \geq \sum_{i \in
C^{\sigma\tau}} \mu(q_i^{\sigma\tau}) \cdot \mu(p_{\sigma (\ell
- 1)}^{\hat\sigma \hat\tau })$$ (cf (6.4) and (6.5)),
and, by symmetry (because our conditions are relatively quadratic),
$$\delta^{-2} \cdot a^{\sigma\tau} \cdot \epsilon \cdot m \cdot z_{m-1}
\geq \sum_{i \in C^{\sigma\tau} } \mu(p_i^{\sigma\tau}) \cdot
\mu(q_{\tau(\ell - 1 )}^{\hat\sigma \hat\tau 
}).$$ Also $\delta \cdot\sum_{i \in j^{\sigma\tau}} 
\mu(p_i^{\sigma\tau}) >
\sum_{i \in A^{\sigma\tau}} \mu (p_i^{\sigma\tau} \cap p)$
and $\delta \cdot\sum_{i\in j^{\sigma\tau}} \mu(q_i^{\sigma\tau}) > \sum_{i \in
B^{\sigma\tau}} \mu (q_i^{\sigma\tau} \cap q)$. So it suffices
to calculate $\sum_{i \in B^{\sigma\tau}, j \in A^{\sigma\tau}} \mu(
p_i^{\sigma\tau} \cap p ,
q_j^{\sigma\tau} \cap q)$. 
\par
For this, we make again backwards induction on $\ell$. Assume
$\ell = m-1$. 
Then the disjointness of $B^{\sigma\tau}$ and $A^{\sigma\tau}$ 
(see (6.1) --- (6.3)) implies that $$\sum_{i \in B^{\sigma\tau} , j \in A^{\sigma\tau}}
\mu(p \cap p_i^{\sigma\tau} , q \cap q_j^{\sigma\tau})
< {1 \over 4} \cdot \mu(p_{\sigma(m-2)}^{\hat\sigma 
\hat\tau }, q_{\tau(m-2)}^{\hat\sigma 
\hat\tau }).$$ Now it follows from the
discussion in the preceding paragraph (and (5.2)) that formula (8.1) holds
for $\ell = m-1$. \par
So assume the claim has been proved for $\ell +1 \leq m-1$.
Let $\sigma\tau$ be nice of length $
\ell$. Then (6.6) $$\sum_{ij , 
\sigma \hat{\;} \langle i \rangle
\tau \hat{\;} \langle j \rangle {\rm not \; nice}} a^{\sigma \hat{\;}
\langle i \rangle \tau \hat{\;} \langle j \rangle} \leq
{ y_\ell \over z_{\ell + 1}} \cdot \mu(p_{\sigma(\ell - 1)}^{\hat\sigma 
\hat\tau }, q_{\tau(\ell - 1)}^{\hat\sigma
\hat\tau}). \eqno (8.2)$$ 
And by induction (and the disjointness of $B^{\sigma\tau}$ and
$A^{\sigma\tau}$) we have $$\eqalignno{\sum_{i \in B^{\sigma\tau} ,
j \in A^{\sigma\tau}, \sigma\hat{\;}\langle i \rangle \tau\hat{\;}
\langle j \rangle {\rm nice}} \mu(p \cap p_i^{\sigma\tau},
& q \cap q_i^{\sigma\tau}) <  ( {1 \over 4^{m - \ell - 1}}
+ \zeta_{\ell + 1}) \cdot \sum_{i \in B^{\sigma\tau},  j \in A^{\sigma\tau}}
\mu(p_i^{\sigma\tau} , q_i^{\sigma\tau})  & \cr & < {1 \over 4} \cdot
({1 \over 4^{m - \ell -1}} + \zeta_{\ell + 1}) \cdot
\mu(p_{\sigma(\ell - 1)}^{\hat\sigma \hat\tau
}, q_{\tau(\ell - 1 )}^{\hat\sigma 
\hat\tau }). & (8.3) \cr}$$
Putting everything ((5.2), the first paragraph
of this proof, (8.2), (8.3)) together we get again that formula (8.1)
holds. $\qed$
\bigskip
The main claim 2.6 now follows from claims 2.7 and 2.8. $\qed$
\bigskip
{\sanse 2.9} {\it Proof of Theorem 3 from the Main Claim 2.6.}
As remarked in 2.5 we let $N := \{
x \in 2^\omega ; \; \exists^\infty m \; (x \restrict n_m
\in Z_m ) \}$. Then $\sum_m {\vert Z_m \vert \over 2^{n_m}}
\leq \sum_m 2^{-m} < \infty$. Hence $N$ is a null set coded in
$V$. We claim that 
\smallskip
\centerline{$\forces_{\BB
\times \BB}  \breve T \cap N \neq \emptyset$.}
\smallskip
\noindent We first note that it suffices to prove
\smallskip
\centerline{$\forces_{\BB \times \BB} 
\forall \ell \in \omega \; \forall s \in \breve T
\; (lh(s) = n_\ell \Rightarrow \exists^\infty m \geq \ell
\; \exists t \; (lh(t) = n_m \land s \subseteq t \land
t \in \breve T \cap Z_m ))$.}
\smallskip
\noindent For if the latter holds then we can recursively
construct an $x \in \breve T [G] \cap N$ in the
generic extension $V [G]$. \par
So assume that there is a $(p,q) \in \BB \times \BB$,
an $\ell \in \omega$ and an $s$ of length $n_\ell$ such
that
\smallskip
\centerline{$(p,q) \forces s \in \breve T \land \forall m
> \ell \; \forall t \; ((lh(t)=n_m \land s \subseteq t)
\Rightarrow t \not\in \breve T \cap Z_m)$;}
\smallskip
\noindent i.e.
\smallskip
\centerline{$(p,q) \forces s \in \breve T \land \forall
m > \ell \; (Z_m^s \cap \breve T = \emptyset)$.}
\smallskip
\noindent Without loss $(p,q)$ is quadratic. Choose $m \geq \ell$
so large that $\mu (p,q) \geq {1 \over 4^m} + \zeta_m$.
Then 
\smallskip
\centerline{$(p,q) \forces s \in \breve T 
\land Z_m^s \cap \breve T = \emptyset$}
\smallskip
\noindent contradicts the main claim 2.6. $\qed$

\vfill\eject

\noindent{\dunhg $\S$ 3 Final remarks}
\Smallskip
{\sanse 3.1} We discuss the relationship between $\CC * \BB$, $\BB * \CC
\cong \BB \times \CC$ and $\BB \times \BB$. Truss [T 2] proved that
$\CC * \BB$ cannot be completely embedded in $\BB * \CC$ by showing
that the former adds a new uncountable subset of $\omega_1$ containing
no old countable subset whereas the latter does not. In fact, he proved
[T 2, Theorem 3.1] that any uncountable subset of $\omega_1$ in 
$V[r]$, where $r$ is random over $V$, contains a countable subset
in $V$; the rest follows from the fact that $\CC$ has a countable
dense subset. It is easy to see that Truss' argument for $\BB$ can
be generalized to $\BB \times \BB$ so that $\CC * \BB$ cannot be 
completely embedded in $\BB \times \BB$ either.
\par
Another argument for showing that $\CC * \BB$ cannot be completely
embedded in $\BB * \CC$ is by remarking that the former produces
two random reals the sum of which is Cohen (namely, let $c$ be
Cohen over $V$ and $r$ random over $V[c]$; then both $r$ and
$c - r$ are random over $V$) whereas the latter does not
(by [JS 2, 2.3], if we force with $\CC$ over $V[r]$, $r$ random
over $V$, then no new real is random over $V$; so the sum of
two random reals must lie in $V[r]$ and cannot be Cohen).
$\BB \times \BB$ also produces two random reals the sum of which
is Cohen (by [Je 2, part I, 5.9], if $r_0$, $r_1$ are the 
random reals added by $\BB \times \BB$, then $r_0 + r_1$ is
Cohen). So $\BB \times \BB$ cannot be completely embedded in
$\BB * \CC$. \par
On the other hand, Pawlikowski (see the last paragraph of $\S$
3 in [Pa]) proved that $\BB * \CC$ can be completely embedded into
any algebra adding both Cohen and random reals; in particular
$\BB * \CC <_c \CC * \BB$ and $\BB * \CC <_c \BB \times \BB$, where
$<_c$ means {\it is complete subalgebra of}. Hence the only
question left open is whether $\BB \times \BB <_c \CC * \BB$
(Question 1).
\bigskip
{\sanse 3.2} We continue with a remark concerning Question 2.
Namely, suppose there are models of $ZFC$ $M \subseteq N$ such
that $N$ contains both a dominating and a random real over $M$,
but does not contain a perfect set of random reals over $M$.
Without loss, $N = M [r] [d]$, where $r$ is random over $M$,
and $d$ is dominating over $M$. By the $\omega^\omega$-bounding
property of random forcing [Je 2, part I, 3.3 (a)], $d$ is
dominating over $M [r]$. Let $c$ be Cohen over $N$. A result
of Truss [T 1, Lemma 6.1] says that $d + c$ is Hechler over $M [r]$.
\par
(Recall that Hechler forcing $\DD$ is defined as follows.
$\DD := \{ (n,f) ; \; n \in \omega \; \land \; f \in \omega^\omega \}$,
$(n,g) \leq (m,f)$ iff $n \geq m$ and $\forall \ell \in \omega
\; (g(\ell) \geq f(\ell))$ and $f \restrict m = g \restrict m$.
$\DD$ generically adds a dominating real.)
\par
By [JS 2, 2.3] there is no new real random over $M$ in $N [c]$,
in particular, there is no perfect set of random reals over $M$
in $N [c]$, thus showing that $\BB * \DD \cong \BB \times \DD$
does not add a perfect set of random reals. Hence Question 2
is equivalent to
\smallskip
{\dunh Question 2'.} {\it Does $\BB \times \DD$ add a perfect
set of random reals?}
\smallskip
\noindent We will now see that for many forcing notions
$\PP$ adding a dominating
real it is true that $\BB \times \PP$ adds a perfect set of random
reals.
\bigskip
{\sanse 3.3} {\capit Proposition.} {\it Let $\MM$ be Mathias forcing.
Then $\BB \times \MM$ adds a perfect set of random reals.}
\smallskip
{\it Remark.} Mathias forcing is defined as follows.
$\MM := \{ (s,S) ; \; s \in \omega^{< \omega} \; \land \;
S \in [\omega]^\omega \; \land \; \max s < \min S \}$,
$(t,T) \leq (s,S)$ iff $t \supseteq s$ and $T \subseteq S$
and $\forall n \in dom(t) \setminus dom(s) \; (t(n) \in S)$.
\smallskip
{\it Sketch of proof.} In $V[G]$, where $G$ is $\BB \times
\MM$-generic over $V$, let $r$ be the random real and $d$
the Mathias real (which is dominating). We claim that
$T: = \{ f \in 2^\omega ; \; \forall n \in \omega \;
(f \restrict [d(n),d(n+1)) = r \restrict [d(n),d(n+1)) \;
\lor \; f \restrict [d(n), d(n+1)) = (1 - r) \restrict
[d(n), d(n+1))) \}$ is a perfect set of reals random over $V$
in $V[G]$. \par
We show that given a null set $N \in V$ and a condition
$(B, (s,S)) \in \BB \times \MM$, there is a $(B' , (s ,S'))
\leq (B,(s,S))$ such that
\smallskip
\centerline{$(B' , (s,S')) \forces \breve T \cap N = \emptyset$,}
\smallskip
\noindent where $\breve T$ is a name for the perfect set defined
above. First note that by [Ba] there are partitions $\{ I_i ; \;
i \in \omega \}$ and $\{ I_i ' ; \; i \in \omega \}$ of
$\omega$ into finite intervals with $\max (I_i) < \min (I_j)$,
$\max (I_i ') < \min (I_j ')$ for $i < j$, sequences 
$\langle J_i ; \; i \in \omega \rangle$ and $\langle J_i ' ;
\; i \in \omega \rangle$ such that $J_i \subseteq 2^{I_i}$,
$J_i ' \subseteq 2^{I_i '}$, $\sum_{i \in \omega} {\vert J_i
\vert  \over 2^{\vert I_i \vert} } < \infty$, $\sum_{i \in
\omega} {\vert J_i ' \vert \over 2^{\vert I_i ' \vert}} <
\infty$ and $N \subseteq \{ f \in 2^\omega ; \; \exists^\infty n
\; ( f \restrict I_n \in J_n ) \} \cup \{ f \in 2^\omega ; \;
\exists^\infty n \; ( f \restrict I_n ' \in J_n ' ) \}$.\par
Now find $S' \subseteq S$ such that for all $n \in \omega$,
$\vert S' \cap I_n \vert \leq 1$ and $\vert S' \cap I_n '
\vert \leq 1$. Let $i_n$ be the unique element of $S' \cap I_n$,
and $i_n '$ the unique element of $S' \cap I_n '$ (if it
exists --- if not, let $i_n \in I_n$ and $i_n ' \in I_n '$
be arbitrary). Set $K_n := \{ s \in 2^{I_n} ; \; s \in J_n \;\lor\;
1-s \in J_n \;\lor\; (s \restrict i_n) \cup (1-s \restrict 
[i_n , \infty) ) \in J_n \;\lor\; (1-s \restrict i_n) \cup
(s \restrict [i_n ,\infty) ) \in J_n \}$; similarly we
define $K_n '$. We choose $B' \leq B$ such that $B' \cap
( \{ f \in 2^\omega ; \; \exists^\infty n \; (f \restrict I_n
\in K_n ) \} \cup \{ f \in 2^\omega ; \; \exists^\infty
n \; (f \restrict I_n ' \in K_n ') \} )= \emptyset$. We leave
it to the reader to verify that this works. $\qed$
\smallskip
We note that a similar argument works for Laver forcing, for
Mathias forcing with a q-point ultrafilter etc. But it is unclear
whether $\BB * \MM$ adds a perfect set of random reals.
\vfill\eject

\noindent{\dunhg References}
\Smallskip
\itemitem{[Ba]} {\capit T. Bartoszy\'nski,} {\it On covering of real
line by null sets,} Pacific Journal of Mathematics, vol. 131 (1988),
pp. 1-12.
\smallskip
\itemitem{[BJ 1]} {\capit T. Bartoszy\'nski and H. Judah,} {\it
Jumping with random reals,} Annals of Pure and Applied Logic,
vol. 48 (1990), pp. 197-213.
\smallskip
\itemitem{[BJ 2]} {\capit T. Bartoszy\'nski and H. Judah,}
{\it Measure and category: the asymmetry,} forthcoming book.
\smallskip
\itemitem{[Je 1]} {\capit T. Jech,} {\it Set theory,} Academic Press,
San Diego, 1978.
\smallskip
\itemitem{[Je 2]} {\capit T. Jech,} {\it Multiple forcing,} 
Cambridge University Press, Cambridge, 1986.
\smallskip
\itemitem{[JS 1]} {\capit H. Judah and S. Shelah,} {\it The 
Kunen-Miller chart (Lebesgue measure, the Baire property,
Laver reals and preservation theorems for forcing),} Journal
of Symbolic Logic, vol. 55 (1990), pp. 909-927.
\smallskip
\itemitem{[JS 2]} {\capit H. Judah and S. Shelah,} {\it Around
random algebra,} Archive for Mathematical Logic, vol. 30 (1990),
pp. 129-138.
\smallskip
\itemitem{[JS 3]} {\capit H. Judah and S. Shelah,} {\it Adding
dominating reals with random algebra,} to appear.
\smallskip
\itemitem{[Ka]} {\capit A. Kamburelis,} {\it Iterations of Boolean
algebras with measure,} Archive for Mathematical Logic, vol. 29
(1989), pp. 21-28.
\smallskip
\itemitem{[Mi]} {\capit A. Miller,} {\it Additivity of measure implies
dominating reals,} Proceedings of the American Mathematical Society,
vol. 91 (1984), pp. 111-117.
\smallskip
\itemitem{[Ox]} {\capit J. C. Oxtoby,} {\it Measure and category,}
Springer, New York Heidelberg Berlin, 2nd edition, 1980.
\smallskip
\itemitem{[Pa]} {\capit J. Pawlikowski,} {\it Why Solovay
real produces Cohen real,} Journal of Symbolic Logic,
vol. 51 (1986), pp. 957-968.
\smallskip
\itemitem{[T 1]} {\capit J. Truss,} {\it Sets having calibre 
$\aleph_1$,} Logic Colloquium 76, North-Holland, Amsterdam,
1977, pp. 595-612.
\smallskip
\itemitem{[T 2]} {\capit J. Truss,} {\it The noncommutativity of
random and generic extensions,} Journal of Symbolic Logic,
vol. 48 (1983), pp. 1008-1012.

\vfill\eject\end